\documentclass[11pt,a4paper,leqno]{article}
\usepackage{a4wide}
\usepackage{latexsym}
\usepackage{amsmath}
\usepackage{amssymb}
\usepackage{stackrel}  
\usepackage{color}
\usepackage{graphicx}
\usepackage[linesnumbered,ruled,vlined]{algorithm2e}

\newtheorem{defin}{Definition}
\newtheorem{lemma}{Lemma}
\newtheorem{prop}{Proposition}
\newtheorem{theo}{Theorem}
\newtheorem{corol}{Corollary}
\newtheorem{example}{Example}
\newenvironment{proof}{\medskip\par\noindent{\bf Proof}}{\hfill $\Box$
\medskip\par}

\newcommand{\C}{\mathbb{C}}
\newcommand{\N}{\mathbb{N}}
\newcommand{\R}{\mathbb{R}}

\begin{document}
\title{Entire solutions of linear systems of moment differential equations and related asymptotic growth at infinity}
\author{{\bf A. Lastra}\\
Universidad de Alcal\'a, Dpto. F\'isica y Matem\'aticas,\\
Alcal\'a de Henares, Madrid, Spain\\
{\tt alberto.lastra@uah.es}}
\date{}
\maketitle
\thispagestyle{empty}
{ \small \begin{center}
{\bf Abstract}
\end{center}

The general entire solution to a linear system of moment differential equations is obtained in terms of a moment kernel function for generalized summability, and the Jordan decomposition of the matrix defining the problem. 

The growth at infinity of any solution of the system is also determined, both globally and also following rays to infinity, determining the order and type of such solutions.

\smallskip

}
\bigskip \bigskip

\section{Introduction}

The main purpose of the present work is to study the explicit description of the general solutions to a linear system of moment differential equations of the form
\begin{equation}\label{e1intro}
\partial_{m}y(z)=Ay(z),
\end{equation}
where $A\in\C^{n\times n}$ is a constant matrix and $y(z)=(y_1(z),\ldots,y_n(z))$ is a vector of unknown functions, for some positive integer $n\ge1$. Here, $m=(m(p))_{p\ge0}$ is a fixed sequence of moments (see Definition~\ref{defi128}), and $\partial_my(z)=(\partial_m(y_1),\ldots,\partial_m(y_n))$ stands for the moment derivative of the vector $y(z)$ (see \ref{e158}), generalizing the classical derivation operator.

The studies of moment differential equations were firstly developed by W. Balser and M. Yoshino in their former work~\cite{bayo}, and have attracted the attention of many researchers in the last years due to the versatility of moment differentiation. Recently, several summability results of the formal solutions to moment partial differential equations have been obtained by S. Michalik~\cite{mi} and by S. Michalik, M. Suwinska and the author~\cite{lamisu}, and also in the final section of the chapter~\cite{sanzrev}, by J. Sanz. Also, some advances have been made on the study of moment integro-differential equations~\cite{lamisu3}. The knowlegdge of the growth of the coefficients of the formal solutions to moment partial differential equations is also of great interest as a first step towards summability results~\cite{su,lamisu2}, and also the Stokes phenomenon~\cite{mitk}.

We also refer to the work~\cite{lamasa2}, where S. Malek, J. Sanz and the author determine summability properties of the solutions of a family of singularly perturbed partial differential equations in the complex domain whose coefficients are $\mathbb{M}$-sums, for some strongly regular sequence $\mathbb{M}$, and also the work~\cite{jikalasa} on multisummability results in Carleman ultraholomorphic classes by means of nonzero proximate orders.

As a matter of fact, versatility of moment differentiation particularizes on the classical differentiation corresponding to the moment sequence $m=(p!)_{p\ge0}$. Moreover, the operator $\partial_{m}$ is quite related to the Caputo $1/s$-fractional differential operator when fixing the sequence of moments $m=(\Gamma(1+\frac{p}{s}))_{p\ge0}$. The theory regarding this concrete choice for the moment sequence is currently being developed due to its applications in many fields of research. Stability properties of systems of fractional differential equations have also been widely studied in literature. In 1996, D. Matignon~\cite{matignon}, describes the stability properties of systems of Caputo fractional differential equations of the form
\begin{equation}\label{e87}
{}^C D^{1/k}_{z}y(z)=Ay(z),
\end{equation}
where ${}^C D^{1/s}_{z}$ stands for Caputo fractional derivative of order $1/s$. Such properties are given in terms of the eigenvalues of the matrix of the system $A$, but also of $s$. The author determines that the solutions of the previous system are asymptotically stable if the eigenvalues $\lambda$ of $A$ satisfy $|\arg(\lambda)|>\frac{\pi}{2s}$. This result is coherent with our asymptotic study of the solutions in Section~\ref{sec5}. This mentioned asymptotic behavior of the solution is due to the appearance of kernel functions associated to summability which construct the analytic solution of the system of moment differential equations. More precisely, we prove (see Theorem~\ref{teodiag} and Theorem~\ref{teonodiag}) that $E(\lambda z)v$ is a solution of (\ref{e1intro}), where $\lambda\in\hbox{spec}(A)$ and $v\in\hbox{ker}(A-\lambda I_n)$ is an associated eigenfunction. $E$ is an entire function, converging to zero in the sector $\{z\in\C:|\arg(z)|>\pi\omega(\mathbb{M})/2\}$, for some positive number $\omega(\mathbb{M})$ related to the moment sequence $m$ (see Definition~\ref{defi128}). We observe both sectors  coincide in the framework of Caputo fractional differential operators.

In practice, the linearized problem (\ref{e87}), much more easier to handle, is usually considered to study the solutions to a nonlinear system, as it is the case of~\cite{ucoz} where the authors study the linearized associated system of the form (\ref{e87}) to study a fractional model of cancer-immune system, providing the asymptotic stability of the solutions from the Jacobian matrix of the initial nonlinear system. We also refer to~\cite{aghror}, where the authors study stability properties of Caputo fractional differential equations under the action of impulses. Other recent results and applications considering linear and semilinear systems of fractional differential equations with constant coefficients are also studied in~\cite{maab,maon}. We also refer to the work by B. Bonilla, M. Rivero and J. J. Trujillo~\cite{bonillaetal} where the authors study systems of linear fractional differential equations with constant coefficients, and the statements of Theorem~\ref{teodiag} and Theorem~\ref{teonodiag} are displayed in terms of the exponential matrix in this framework. 

Given a strongly regular sequence $\mathbb{M}$ admitting a kernel function $E$ for $\mathbb{M}$-summability (see Definition~\ref{defi132} and Definition~\ref{defi128}), the properties associated to the entire function $E$ allow us to construct the general solution of (\ref{e1intro}) in Theorem~\ref{teodiag} (for a diagonal matrix $A$) and Theorem~\ref{teonodiag} (under general settings). This general solution, depending on $n$ arbitrary constants, is obtained by appropriate manipulation of $E$, the eigenvalues and generalized eigenvectors of $A$, together with the Jordan canonical form associated to $A$. Once obtained the explicit entire solutions of the system (\ref{e1intro}), we study their asymptotic behavior at infinity with the help of the classical theory of growth of entire functions. This leads us to the knowledge of the global growth of the entire solutions of (\ref{e1intro}) at infinity not only in a global sense (Theorem~\ref{teo5}), but also determining the radial growth of the solutions at infinity (Theorem~\ref{teo6}). The main novelty of the present work is twofold. On the one hand, a closed explicit expression for the entire solutions to linear systems of moment partial differential equations of the form (\ref{e1intro}) is obtained (Theorem~\ref{teodiag} and Theorem~\ref{teonodiag}) in the general framework of moments related to strongly regular sequences $\mathbb{M}$ admitting kernel functions for $\mathbb{M}$-summability (see Section~\ref{sec3}). On the other hand, we also determine the order and type of growth associated to such entire solutions at infinity both globally (Theorem~\ref{teo5}) and also following rays approaching to infinity (Theorem~\ref{teo6}).

The paper is organized as follows: After fixing some notation considered in the work (Section~\ref{sec2}) we recall the main facts on strongly regular sequences and generalized summability (Section~\ref{sec3}). The general explicit entire solutions to (\ref{e1intro}) are constructed in Section~\ref{sec4}, and their asymptotic behavior at infinity is obtained in Section~\ref{sec5}, with a special attention to the radial growth at infinity of the entire solutions (Section~\ref{sec51}).

\section{Notation}\label{sec2}

$\mathcal{R}$ denotes the Riemann surface of the logarithm. Let $d\in\R$, and $\theta>0$. We write $S_d(\theta)$ for the set
$$S_d(\theta):=\left\{z\in\mathcal{R}:|\arg(z)-d|<\frac{\theta}{2}\right\}.$$
$I_n$ stands for the identity operator in $\C^{n\times n}$ for some $n\in\N:=\{1,2,\ldots\}$. Given a set $C$, we write $\#(C)$ for the cardinal number of $C$.

\noindent $D(z,r)$ stands for the open disc centered at $z\in\C$ and radius $r>0$.

\noindent $\mathbb{C}[[z]]$ denotes the set of formal power series in $z$ with complex coefficients. Given an open set $U\subseteq \C$, we write $\mathcal{O}(U)$ for the set of holomorphic functions on $U$.

\section{Strongly regular sequences and generalized summability}\label{sec3}

This section is devoted to recall the main facts on the elements involved in the results of the present work with respect to the theory of generalized summability. We mainly focus on the definition and main properties of a strongly regular sequence $\mathbb{M}$, and the concept of kernel functions for $\mathbb{M}$-summability, leading to the notion of sequence of moments, and other tools associated to generalized summability and generalized differential operators. The main problem under study will be written in terms of these tools and notions.

\subsection{Strongly regular sequences}

The notion of strongly regular sequence is due to V. Thilliez~\cite{thilliez}, generalizing Gevrey sequences.

\begin{defin}\label{defi132}
Let $\mathbb{M}:=(M_p)_{p\ge0}$ be a sequence of positive real numbers such that $M_0=1$ and satisfies the following properties:
\begin{itemize}
\item[(lc)] $\mathbb{M}$ is logarithmically convex, i.e. $M_{p}^2\le M_{p-1}M_{p+1}$ for every $p\ge1$.
\item[(mg)] $\mathbb{M}$ is of moderate growth, i.e. there exists $A_1>0$ with $M_{p+q}\le A_1^{p+q}M_pM_q$, for every pair of integers $p,q\ge0$.
\item[(snq)] $\mathbb{M}$ is non-quasianalytic, i.e. there exists $A_2>0$ such that
$$\sum_{q\ge p}\frac{M_q}{(q+1)M_{q+1}}\le A_2\frac{M_p}{M_{p+1}}$$
for every $p\ge0$.
\end{itemize}
Any sequence $\mathbb{M}$ satisfying the previous properties is known as a strongly regular sequence.
\end{defin}

Examples of strongly regular sequences can be found in the literature in many applications. For example Gevrey sequences, defined by $\mathbb{M}_{\alpha}:=(p!^{\alpha})_{p\ge0}$ for some $\alpha>0$, are of great importance in the theory of summability of formal solutions to ordinary and partial differential equations (see~\cite{ba2}, and the references therein). Given $\alpha>0$ and $\beta\in\R$, the sequence $\mathbb{M}_{\alpha,\beta}:=(p!^{\alpha}\prod_{m=0}^{p}\log^{\beta}(e+m))_{p\ge0}$ is also a strongly regular sequence, the first terms being slightly modified if $\beta<0$. This modification does not interfere with the spaces of functions related to that sequence. As a matter of fact, the particular case $\alpha=-\beta=1$ is linked to the $1+$ level in the study of difference equations~\cite{immink1,immink2}.

Associated to a strongly regular sequence $\mathbb{M}$ one can define the non-decreasing continuous function $M:[0,\infty)\to[0,\infty)$ by $M(0)=0$ and 
\begin{equation}\label{e94}
M(t):=\sup_{p\ge0}\log\left(\frac{t^p}{M_p}\right),\quad t>0
\end{equation}
(see~\cite{mandel}). In the work~\cite{sanz}, J. Sanz proves that the order of $M$, previously used in the theory of meromorphic functions in~\cite{goos}, and defined by
$$\rho(M):=\limsup_{r\to\infty}\max\left\{0,\frac{\log(M(r))}{\log(r)}\right\},$$
turns out to be a positive real number which provides the inverse of the limit opening of injectivity of the asymptotic Borel map. Namely, the existence of nonzero functions asymptotically flat at the vertex of a sectorial region is guaranteed whenever the opening of the sectorial region is at most $\pi\omega(\mathbb{M})$, with $\omega(\mathbb{M}):=1/\rho(M)$ (see Corollary 3.16 in~\cite{jisash}). 


\subsection{Generalized summability}

The following definitions and properties correspond to the theory of generalized summability, associated to a strongly regular sequence, extending the classical notion of Gevrey summability and the moment summability methods developped by W. Balser, which can be found in~\cite{ba2}, section 5.5. The generalization of that theory to strongly regular sequences was put forward by J. Sanz in~\cite{sanz}. We recall in this section some of the main concepts for the sake of completeness, which will be used in the sequel.

\begin{defin}\label{defi128}
Let $\mathbb{M}$ be a strongly regular sequence, with $0<\omega(\mathbb{M})<2$, and define the function $M$ according to (\ref{e94}). We say that $\mathbb{M}$ admits a pair of kernel functions for $\mathbb{M}$-summability if there exist two functions, $e$ and $E$ satisfying the following conditions:
\begin{enumerate}
\item $e$ is holomorphic on the sector $S_0(\omega(\mathbb{M})\pi)$, and $e(z)/z$ is locally uniformly integrable at the origin meaning there exists $t_0>0$, and for every $z_0\in S_0(\omega(\mathbb{M})\pi)$ there exists $r_0=r_0(z_0)>0$ such that $D(z_0,r_0)\subseteq S_0(\omega(\mathbb{M})\pi)$ and 
$$\int_0^{t_0}\sup_{z\in D(z_0,r_0)|}e(t/z)|\frac{dt}{t}<\infty.$$
In addition to this, for every fixed $\delta>0$, there exist $k_1,k_2>0$ with
$$|e(z)|\le k_1\exp\left(-M\left(\frac{|z|}{k_2}\right)\right), \qquad z\in S_0(\omega(\mathbb{M})\pi-\delta).$$
\item $e(|z|)\in\R$ for every $z\in S_0(\omega(\mathbb{M})\pi)$.
\item $E$ is an entire function, and there exist $k_3,k_4>0$ with
$$|E(z)|\le k_3\exp\left(M\left(\frac{|z|}{k_4}\right)\right),\quad z\in\C.$$
Moreover, there exists $\beta>0$ such that for every $0<\theta<2\pi-\omega(\mathbb{M})\pi$ and $R>0$, there exists $k_5>0$ with
\begin{equation}\label{e145}
|E(z)|\le k_5|z|^{-\beta},\quad z\in S_{\pi}(\theta)\hbox{ with }|z|\ge R.
\end{equation}
\item $E$ admits the representation
\begin{equation}\label{e143}
E(z)=\sum_{p\ge0}\frac{z^p}{m(p)},\quad z\in\C,
\end{equation}
where $m$ is the moment function associated with $e$
\begin{equation}\label{e151}
m(z):=\int_0^{\infty}t^{z-1}e(t)dt,
\end{equation}
for every $z\in\{\omega\in\C:\hbox{Re}(\omega)\ge 0\}$. The function $m$ turns out to be continuous in $\{z\in\C:\hbox{Re}(z)\ge 0\}$ and holomorphic in $\{z\in\C:\hbox{Re}(z)> 0\}$. The sequence $(m(p))_{p\ge0}$ is known as the sequence of moments associated with the pair of kernel functions.
\end{enumerate}
\end{defin}

The existence of a pair of kernel functions for $\mathbb{M}$-summability is not always guaranteed. However, it is when departing from a strongly regular sequence admitting a nonzero proximate order in the sense of E. Lindel\"of (see~\cite{jisash2,lamasa}), which is the case of sequences appearing in applications, as stated in~\cite{sanz}. We also refer to~\cite{jisash2}, where the authors determine conditions for a strongly regular sequence for admitting a nonzero proximate order. In particular, every sequence of the form $\mathbb{M}_{\alpha,\beta}$ for any $\alpha>0$ and $\beta\in\R$ as defined above, and also Gevrey sequence $\mathbb{M}_{\alpha}$ for every $\alpha>0$ admit  nonzero proximate order (see Theorem 3.6 and Example 3.7 in~\cite{jisash2}) In addition to this, $\omega(\mathbb{M}_{\alpha,\beta})=\alpha$ and $\omega(\mathbb{M}_{\alpha})=\alpha$.

\begin{defin}
A nonzero proximate order $\rho(t)$ is a nonnegative continuously differentiable function defined in an interval of the form $(c,\infty)$ for some $c\in\R$ such that
\begin{equation}\label{e159}
\lim_{r\to\infty}\rho(r)=\rho,
\end{equation}
for some $\rho\in\R$, and 
\begin{equation}\label{e169}
\lim_{r\to\infty}r\rho'(r)\ln(r)=0.
\end{equation}
\end{defin}

In the case of $\mathbb{M}_{\alpha,\beta}$ and $\mathbb{M}_{\alpha}$ with $\alpha>0$ and $\beta\in\R$, one has $\rho=1/\alpha$.

The following result can be found in~\cite{sanz}, Proposition 5.8.
\begin{prop}\label{propequiv}
Given a strongly regular sequence $\mathbb{M}=(M_p)_{p\ge0}$. Assume that $\mathbb{M}$ admits a pair of kernel functions for $\mathbb{M}$-summability, and let $m=(m(p))_{p\ge0}$ be the associated moment sequence. Then, $\mathbb{M}$ and $m$ are equivalent sequences, i.e. there exist $c_1,c_2>0$ such that 
$$c_1^pM_p\le m(p)\le c_2^pM_p,\qquad p\ge0.$$
\end{prop}

Given a sequence of positive real numbers, or in practice a sequence of moments $m=(m(p))_{p\ge0}$, Balser and Yoshino~\cite{bayo} defined the moment derivative operator $\partial_{m,z}:\mathbb{C}[[z]]\to\mathbb{C}[[z]]$ by
\begin{equation}\label{e158}
\partial_{m,z}\left(\sum_{p\ge0}\frac{a_p}{m(p)}z^p\right)=\sum_{p\ge0}\frac{a_{p+1}}{m(p)}z^p.
\end{equation}
Its definition can be naturally extended to any holomorphic function defined on some neighborhood of the origin when considering the Taylor expansion of the function. Such formal operator corresponds to the usual derivation for the sequence $m=(p!)_{p\ge0}$ and it is quite related to the Caputo $1/s$-fractional differential operator ${}^C D^{1/k}_{z}$ when considering the sequence of moments $m=(\Gamma(1+\frac{p}{k}))_{p\ge0}$. Indeed, for the previous moment sequence one has that
\begin{equation}\label{e155}
( \partial_{m,z}f)(z^{1/k})={}^C D^{1/k}_{z}(f(z^{1/k}))
\end{equation}
for every $f\in\mathbb{C}[[z]]$. We refer to~\cite{michalik12} (see Definition 5 and Remark 1) for further details. In addition to this, the definition of moment differential operators can be extended not only to holomorphic functions near the origin, but also to summable functions. More precisely, given a strongly regular sequence $\mathbb{M}$ admitting a nonzero proximate order, and a function $f$ defined on a finite sector $S$ of opening larger than $\pi\omega(\mathbb{M})$ and bisecting direction $d\in\R$, which admits the formal power series $\hat{f}\in\mathbb{C}[[z]]$ as its $\mathbb{M}$-asymptotic expansion in $S$ (observe from Corollary 4.12~\cite{sanz} this function is unique under this property), then one can define $\partial_{m,z}f$ as the unique function admitting $\partial_{m,z}(\hat{f})$ as its $\mathbb{M}$-asymptotic expansion in a finite sector of opening larger than $\pi\omega(\mathbb{M})$ and bisecting direction $d\in\R$. We refer to~\cite{lamisu} for a more detailed description of this result.

\section{Entire solutions to a linear system of moment differential equations}\label{sec4}

In this section, we state the first main result of the present work. Namely, we give the precise analytic expression of the solution to certain linear system of moment differential equations. The asymptotic behavior of such solutions at infinity will be analyzed in a subsequent section.

Let $\mathbb{M}$ be a strongly regular sequence admitting a nonzero proximate order. This entails the existence of an associated pair of kernel functions for $\mathbb{M}$-summability, say $e$ and $E$, satisfying the properties enumerated in Definition~\ref{defi128}. Let $m=(m(p))_{p\ge0}$ be the sequence of moments associated to the previous kernel functions. Let $n\ge 1$, and $A\in\mathbb{C}^{n\times n}$ be a $n\times n$ matrix with complex coefficients. We consider the moment differential equation
\begin{equation}\label{epralhom}
\partial_my=Ay,
\end{equation}
where $y=(y_1,\ldots,y_n)^{T}$ is a vector of unknown functions, with $y_j=y_j(x)$ for $1\le j\le n$, and $\partial_{m}y$ stands for the vector $(\partial_{m}y_1,\ldots, \partial_{m}y_n)^{T}$.

\begin{lemma}\label{lema1}
Let $(z_0,y_0)\in\C^{1+n}$. The Cauchy problem 
\begin{equation}\label{e176}
\left\{ \begin{array}{lcc}
             \partial_my&=Ay\\
              y(z_0)&=y_0
             \end{array}
   \right.
\end{equation}
admits a unique solution which is and entire function.
\end{lemma}
\begin{proof}
One can set $z_0=0$ without loss of generality. Indeed, if $y$ is a solution of the equation with Cauchy data located at the origin, then $y(z-z_0)$ is a solution of the same equation with its Cauchy data located at $z=z_0$.

We first show that (\ref{e176}) admits a unique formal power series solution. Let us write $A=(a_{jk})_{1\le j,k,\le n}$, and assume that $y(z)=\sum_{p\ge0}\tilde{y}_p\frac{z^p}{m(p)}\in\C^{n}[[z]]$ is the formal solution of (\ref{e176}). We write $\tilde{y}_p=(\tilde{y}_{p,1},\ldots,\tilde{y}_{p,n})$ for every $p\ge0$. Then, one has that $\tilde{y}_{p+1,j}=\sum_{k=1}^{n}a_{jk}\tilde{y}_{p,j}$ for every $1\le j\le n$ and $p\ge0$. Let $\alpha:=\max_{1\le j,k\le n}|a_{j,k}|$, and put $\left\|\tilde{y}_p\right\|:=\sum_{j=1}^{n}|\tilde{y}_{p,j}|$ for every $p\ge0$. For all $p\ge1$ and $1\le j\le n$ one has
$$\left\|\tilde{y}_{p+1}\right\|=\sum_{j=1}^{n}|\tilde{y}_{p+1,j}|\le \sum_{j=1}^{n}\sum_{k=1}^{n}|a_{jk}||\tilde{y}_{p,j}|\le \alpha n\left\|\tilde{y}_{p}\right\|.$$
We define the sequence $(c_p)_{p\ge0}$ by $c_0:=\left\|\tilde{y}_0\right\|$, and for every $p\ge0$ we put $c_{p+1}:= \alpha n c_p$. This entails that $c_p=(\alpha n)^{p}c_0$ for every $p\ge0$. In addition to this, it is straight to check that $\left\|\tilde{y}_p\right\|\le c_p$ for every $p\ge0$. Therefore, one has that $|\tilde{y}_{p,j}|\le (\alpha n)^{p}\left\|\tilde{y}_0\right\|$ for all $1\le j\le n$. The series $\sum_{p\ge0}(\alpha n |z|)^p\frac{1}{m(p)}$ has infinite radius of convergence. Indeed, it coincides with $E(\alpha n|z|)$ (see (\ref{e143})). Therefore, the formal power series $y(z)$ represents an entire function.
\end{proof}

A proof of the following result follows analogous arguments as in the case of the usual derivation. We only sketch its proof for the sake of completeness.

\begin{lemma}
The set of solutions to (\ref{epralhom}) is a subspace of $(\mathcal{O}(\C))^{n}$ of dimension $n$.
\end{lemma}

\begin{proof}
 Let $y$ be a solution to (\ref{epralhom}), which is defined for some $z_{0}\in\C$. In view of Lemma~\ref{lema1}, one derives that $y$ belongs to $(\mathcal{O}(\C))^n$. It is straight to check that any linear combination of two solutions to (\ref{epralhom}) is also a solution to (\ref{epralhom}). Let $\{\omega_1,\ldots,\omega_n\}$ be a basis of $\C^{n}$ and consider the Cauchy problem
\begin{equation}\label{e203}
\left\{ \begin{array}{lcc}
             \partial_my&=Ay\\
              y(z_0)&=\omega_j
             \end{array}
   \right.
\end{equation}
for every $1\le j\le n$. Lemma~\ref{lema1} guarantees the existence of a unique entire solution to (\ref{e203}), say $y_j$, for each $1\le j\le n$ . The set $\{y_1,\ldots,y_n\}$ determines a basis of the vector space of solutions to (\ref{epralhom}).
\end{proof}

The properties of the kernel functions described in Definition~\ref{defi128} allow us to obtain an explicit formula for the general solution of (\ref{epralhom}). First, we find the explicit solutions to (\ref{epralhom}).

\begin{lemma}\label{lema219}
Let $\lambda\in\C$ be an eigenvalue of $A\in\C^{n\times n}$ with associated eigenvector $v\in\C^n$. The following properties hold:
\begin{itemize}
\item[-] If $\lambda\neq 0$, then the function $y(z)=E(\lambda z)v$ is an entire solution of (\ref{epralhom}). 
\item[-] If $\lambda=0$, then the constant function $y(z)=v$ is a solution of (\ref{epralhom}).
\end{itemize} 
\end{lemma}
\begin{proof}
First, assume that $\lambda\neq 0$. The holomorphy properties of the kernel function $E$ guarantee that $y(z)=E(\lambda z)v$ belongs to $(\mathcal{O}(\C))^n$. 
It is straight to check that $\partial_m(E(\lambda z))=\lambda E(\lambda z)$. Indeed,
\begin{equation}\label{e222}
\partial_m(E(\lambda z))=\partial_{m}\left(\sum_{p\ge0}\frac{(\lambda z)^p}{m(p)}\right)=\sum_{p \ge0}\frac{\lambda^{p+1}z^p}{m(p)}=\lambda E(\lambda z).
\end{equation}
In view of (\ref{e222}) one has that
$$A y(z)=A(E(\lambda z) v)=E(\lambda z) (A v)=E(\lambda z)\lambda v=\partial_m(E(\lambda z))v=\partial_m(E(\lambda z)v)=\partial_m(y(z)).$$
This concludes the first part of the  proof. Now, assume that $\lambda=0$. Then, one has that
$$Av=0=\partial_{m}(v).$$
Therefore, the constant vector $y(z)=v$ provides a solution of (\ref{epralhom}).
\end{proof}

\begin{theo}\label{teodiag}
Let $A\in\C^{n\times n}$ be a diagonalizable matrix. Let $\{\lambda_j\}_{1\le j\le k}$, for some $1\le k\le n$ be the set of eigenvalues of $A$, and let $\{v_{j,1},\ldots, v_{j,\ell_j}\}$ be a basis of $\hbox{Ker}(A-\lambda_jI_n)$ for every $1\le j\le k$, and some $\ell_j\ge 1$. Then, the general solution of (\ref{epralhom}) is given by
\begin{equation}\label{e235}
y(z)=\sum_{j=1}^{k}\sum_{p=1}^{\ell_j}C_{j,p}E(\lambda_j z)v_{j,p},
\end{equation}
with $C_{j,p}$ being arbitrary constants.
\end{theo}
\begin{proof}
In view of Lemma~\ref{lema219}, one obtains that for every $1\le j\le k$ and all $1\le p\le \ell_j$ the function $y_{j,p}(z)=E(\lambda_j)v_{j,p}$ is an entire solution of (\ref{epralhom}). The fact that $A$ is diagonalizable guarantees that $\#\{v_{j,p}:1\le j\le k, 1\le p\le \ell_j\}=n$. It only rests to prove that the set $\{E(\lambda_j z)v_{j,p}: 1\le j\le k,\lambda_j\neq0,1\le p\le \ell_j\}$ is linearly independent. This last statement is straight to be checked by taking into account that $\{v_{j,p}:1\le j\le k, 1\le p\le \ell_j\}$ is a basis of $\C^n$. 
\end{proof}

\noindent\textbf{Remark:} Observe that the terms of the sum in (\ref{e235}) which correspond to the eigenvalue $\lambda=0$ are such that $E(\lambda z)\equiv 1$. This will have consequences regarding the asymptotic behavior of the solution at infinity.

The general case makes use of Jordan canonical form. For this purpose, we define the following formal power series.

\begin{defin}
Let $\lambda\in\C$. For every $h\ge0$ we define the formal power series
\begin{equation}\label{e258}
\Delta_hE(\lambda z)=\sum_{p\ge h}{p\choose h}\frac{\lambda^{p-h}z^p}{m(p)}.
\end{equation}
\end{defin}

\begin{lemma}\label{lema288}
The following statements are direct consequence of the definition of the operator $\Delta_{h}$:
\begin{enumerate}
\item  $\Delta_0E(\lambda z)=E(\lambda z)$.
\item $\Delta_1E(\lambda z)=\partial_{\lambda}(E(\lambda z))$.
\item If $\lambda=0$, $\Delta_hE(\lambda z)=z^{h}/m(h)$. 
\item For every $h\ge1$ it holds that
\begin{equation}\label{e271}
(\partial_m-\lambda)(\Delta_hE(\lambda z))=\Delta_{h-1}E(\lambda z).
\end{equation}
\end{enumerate}
\end{lemma}
\begin{proof}
We only give detail of the proof of the last statement, being the other ones straightforward properties derived from the definition of the operator $\Delta_h$. Let $h\ge1$. In view of (\ref{e158}), it holds that
\begin{equation}\label{e280}
\partial_{m}\left(\sum_{p\ge h}{p\choose h}\frac{\lambda^{p-h}z^p}{m(p)}\right)=\sum_{p\ge h}{p \choose h}\frac{\lambda^{p-h}z^{p-1}}{m(p-1)}.
\end{equation}
On the other hand, one has
\begin{multline}
\lambda\Delta_h E(\lambda z)+\Delta_{h-1}E(\lambda z)=\lambda\sum_{p\ge h}{p\choose h} \frac{\lambda^{p-h}z^p}{m(p)}+\sum_{p\ge h-1}{p\choose h-1} \frac{\lambda^{p-(h-1)}z^p}{m(p)}\\
=\sum_{p\ge h}\left[{p\choose h}+{ p \choose h-1}\right] \frac{\lambda^{p-h+1}z^p}{m(p)}+\frac{z^{h-1}}{m(h-1)}\\
=\sum_{p\ge h}{p+1\choose h} \frac{\lambda^{p-h+1}z^p}{m(p)}+\frac{z^{h-1}}{m(h-1)}=\sum_{p\ge h}{p \choose h}\frac{\lambda^{p-h}z^{p-1}}{m(p-1)},
\end{multline}
leading to the conclusion.
\end{proof}

\begin{lemma}
For every $h\ge0$, and every $\lambda\in\C$, the formal power series $\Delta_hE(\lambda z)$ is an entire function. 
\end{lemma}
\begin{proof}
Taking into account statements 1. and 3. of Lemma~\ref{lema288}, the result is clear for $\lambda=0$ and also for $h=0$. Assume that $\lambda\in\C^{\star}$ and let $h\ge 1$. The statement is clear provided that
$$\Delta_hE(\lambda z)=\left[\frac{1}{\lambda^{h}h!}\omega^{h}\left(\frac{d}{dw}\right)^{h}E(\omega)\right]_{w=\lambda z}.$$
\end{proof}

We now describe the general solution of (\ref{epralhom}) in the case of a non-diagonalizable matrix $A$. Let $p_{A}(z)$ be the characteristic polynomial of $A$, i.e. the polynomial $\hbox{det}(A-z I_n)$. 

\begin{theo}\label{teonodiag}
Let $\{\lambda_j\}_{1\le j\le k}$ for some $1\le k\le n$ be the set of eigenvalues of $A$. Assume that $\lambda_j$ is an eigenvalue of algebraic multiplicity $m_j\ge1$, for every $1\le j\le k$. Then, the general solution of (\ref{epralhom}) can be written in the form
\begin{equation}\label{e283}
y(z)=\sum_{j=1}^{k}\sum_{p=1}^{m_j}C_{j,p}y_{j,p}(z),
\end{equation}
where $C_{j,p}$ are arbitrary constants, and $y_{j,p}$ are entire functions determined from the Jordan decomposition of $A$.
\end{theo}



\begin{proof}
Let $1\le j\le k$ and assume that the algebraic multiplicity of $\lambda_j$, say $m_j$, is larger than its geometric multiplicity, say $\ell_j\ge1$. Otherwise, one can proceed as in Theorem~\ref{teodiag}. The functions $y_{j,1},\ldots,y_{j,m_j}$ are constructed as follows. We write $\lambda:=\lambda_j$ for simplicity. Let $\{v_{j,1},\ldots,v_{j,\ell_j}\}$ be a basis of $\hbox{Ker}(A-\lambda_jI_n)$. For every $1\le p\le \ell_j$ we proceed with an analogous construction as that of a classical Jordan block related to $\lambda_j$, see~\cite{shafarevich}, as an example of reference in this direction. We illustrate the procedure for the sake of completeness. Let us write $\tilde{v}_1:=v_{j,p}$. We put
$$y_{j,1}:=E(\lambda z)\tilde{v}_1.$$
We choose $\tilde{u}_2\in\C^n$ such that 
\begin{equation}\label{e292}
A\tilde{u}_2=\lambda \tilde{u}_2+\tilde{v}_1,
\end{equation}
i.e. $\tilde{u}_2\in(A-\lambda I_n)^{-1}(\tilde{v}_1)$. Observe there exists such a vector due to the algebraic and geometric multiplicities of $\lambda$ do not match. Therefore, $\tilde{u}_2\in\ker((A-\lambda I_n)^2)$. We define
$$y_{j,2}:=\Delta_1E(\lambda z)\tilde{v}_1+E(\lambda z)\tilde{u}_2.$$
If the procedure for the Jordan blocking is not concluded, we proceed analogously by choosing $\tilde{u}_3\in (A-\lambda I_n)^{-1}(\tilde{u}_2)$ (therefore $\tilde{u}_3\in\ker((A-\lambda I_n)^3)$), and defining
$$y_{j,3}:=\Delta_2E(\lambda z)\tilde{v}_1+\Delta_1E(\lambda z)\tilde{u}_2+E(\lambda z)\tilde{u}_3.$$
The recursion concludes with
$$y_{j,h_{p,j}}:=\Delta_{h_{p,j}-1}E(\lambda z)\tilde{v}_1+\Delta_{h_{p,j}-2}E(\lambda z)\tilde{u}_2+\ldots+ E(\lambda z)\tilde{u}_{h_{p,j}},$$
where $\Delta_{h}$ is the operator defined in (\ref{e258}), and $h_{p,j}\ge2$ is determined by the Jordan decomposition associated to $A$, and satisfying $\sum_{j=1}^{k}\sum_{p=1}^{\ell_j}h_{p,j}=n$.

First, observe that
$$\partial_{m}(y_{j,1})=\partial_{m}(E(\lambda z))\tilde{v}_1=\lambda E(\lambda z) \tilde{v}_1= E(\lambda z)A\tilde{v}_1=Ay_{j,1}.$$
Also, in view of the property (\ref{e271}) and (\ref{e292}), one has that 
\begin{multline}
\partial_m(y_{j,2})=\partial_m(\Delta_1E(\lambda z)\tilde{v}_1+E(\lambda z)\tilde{u}_2)=\partial_{m}(\Delta_1E(\lambda z)\tilde{v}_1)+\partial_m(E(\lambda z)\tilde{u}_2)\\
=\lambda\Delta_1E(\lambda z)\tilde{v}_1+ E(\lambda z)\tilde{v}_1+\lambda E(\lambda z)\tilde{u}_2=\Delta_{1}E(\lambda z)A\tilde{v}_1+E(\lambda z)A\tilde{u}_2=Ay_{j,2}.
\end{multline}
This entails that $y_{j,2}$ solves (\ref{epralhom}). One can proceed recursively to obtain for all $1\le q\le h_{p,j}$ that $y_{j,q}$ is a solution of (\ref{epralhom}). Indeed, one has
\begin{multline}
\partial_m(y_{j,q})=\partial_m(\Delta_{q-1}E( \lambda z)\tilde{v}_1+\Delta_{q-2}E(\lambda z)\tilde{u}_2+\ldots+E(\lambda z)\tilde{u}_{q})\\
=\partial_m(\Delta_{q-1}E( \lambda z)\tilde{v}_1)+\partial_m(\Delta_{q-2}E(\lambda z)\tilde{u}_2)+\ldots+\partial_m(E(\lambda z)\tilde{u}_{q})\\
=\lambda\Delta_{q-1}E( \lambda z)\tilde{v}_1+\Delta_{q-2}E(\lambda z)\tilde{v}_1+
\lambda\Delta_{q-2}E( \lambda z)\tilde{u}_2+\Delta_{q-3}E(\lambda z)\tilde{u}_2\\
\hfill+\ldots+\lambda\Delta_1E(\lambda z)\tilde{u}_{q-1}+E(\lambda z)\tilde{u}_{q-1}+\lambda E(\lambda z)\tilde{u}_q\\
=\lambda\Delta_{q-1}E(\lambda z)\tilde{v}_1+\Delta_{q-2}E( \lambda z)(\tilde{v}_1+\lambda \tilde{u}_2)+\Delta_{q-3}E( \lambda z)(\tilde{u}_2+\lambda \tilde{u}_3)\hfill\\
\hfill+\ldots+\Delta_1E( \lambda z)(\tilde{u}_{q-2}+\lambda \tilde{u}_{q-1})+E(\lambda z)(\tilde{u}_{q-1}+\lambda \tilde{u}_q)\\
=\Delta_{q-1}E(\lambda z)A\tilde{v}_1+\Delta_{q-2}E( \lambda z)A \tilde{u}_2+\Delta_{q-3}E( \lambda z)A\tilde{u}_3+\ldots+\Delta_1E( \lambda z)A\tilde{u}_{q-1}+E(\lambda z)A\tilde{u}_{q}=Ay_{j,q}.
\end{multline}

It is straight to check that the set $\{y_{j,p}:1\le j\le k,1\le p\le m_j\}$ determines a basis of the vector space of solutions to (\ref{epralhom}) as it is conformed with $n$ vectors which are linearly independent. Indeed, for every $1\le j\le k$ and $1\le p\le \ell_j$ we observe that $\{y_{j,1},\ldots,y_{j,h_{p,j}}\}$ are obtained by a triangular linear combination of the basis of the Jordan block $\{\tilde{v}_1,\ldots, \tilde{u}_{h_{p,j}}\}$, where the terms in the diagonal are given by $E(\lambda z)$. This concludes the proof. 
\end{proof}

In view of Theorem~\ref{teodiag} and Theorem~\ref{teonodiag} one can state an algorithm for the construction of the general solution of (\ref{epralhom}), by means of the Jordan decomposition of the matrix $A$ and the functions $\Delta_{h}E(\lambda z)$ defined in (\ref{e258}). We illustrate the previous theory in several examples.

\begin{example}
Let $\mathbb{M}_{1}=(p!)_{p\ge0}$ be Gevrey sequence of order 1. $\mathbb{M}_1$ is a strongly regular sequence which admits a nonzero proximate order. The associated kernel function $E(z)=\exp(z)$ and for every $h\ge0$ one gets that $\Delta_{h}E(\lambda z)=\frac{z^h}{h!}\exp(\lambda z)$. The results of Theorem~\ref{teodiag} and Theorem~\ref{teonodiag} coincide with the classical theory of solutions to linear systems of differential equations.
\end{example}

\begin{example}\label{ex2}
This theory can be applied to fractional differential equations, where the Caputo $1/s-$ fractional differential operators ${}^C D^{1/k}_{z}$ are involved in the equation. Indeed, the expression (\ref{e155}) relates both problems, that in terms of Caputo $1/s$-fractional derivatives and the moment differential equations, with moment sequence $m=(\Gamma(1+\frac{p}{k}))_{p\ge0}$. Indeed, one observes in this case that 
$$h!\Delta_h E(\lambda t^{1/k})=t^{h/k}\left(\frac{d^{h}}{dz^{h}}E_{1/k}\right)(\lambda t^{1/k}),$$
for all $h\ge0$, where $E_{1/k}(z)$ is Mittag-Leffler function $E_{1/k}(z)=\sum_{p\ge0}\frac{z^p}{\Gamma(1+p/k)}$.

We refer to the work by B. Bonilla, M. Rivero and J. J. Trujillo~\cite{bonillaetal} where the authors study systems of linear fractional differential equations with constant coefficients, and the statements of Theorem~\ref{teodiag} and Theorem~\ref{teonodiag} are displayed in terms of the exponential matrix in this framework.

\end{example}

\section{Asymptotic study of the solutions}\label{sec5}

In this section, we analyze the asymptotic behavior at infinity of the entire solutions of the problem (\ref{epralhom}), giving rise to stability results of such systems of moment differential equations. We maintain the same assumptions as in the previous section for the elements involved in the system (\ref{epralhom}).

First, we study the global growth of the entire solutions at infinity in terms of the eigenvalues of the matrix $A$ in (\ref{epralhom}). For this purpose, we make use of Proposition 4.5,~\cite{komatsu}, which remains valid under more general arguments.

\begin{prop}\label{propk}
Let $\mathbb{M}$ be a strongly regular sequence, and let $f(z)=\sum_{p\ge0}a_pz^p$ be an entire function. Then, the following statements are equivalent:
\begin{itemize}
\item[-] There exist $C_1,C_2>0$ such that $|f(z)|\le C_1\exp(M(C_2|z|))$, for $z\in\C$.
\item[-] There exist $D_1,D_2>0$ such that $|a_p|\le D_1D_2^p/M_p$, for all $p\ge0$.
\end{itemize}
\end{prop}

\begin{prop}\label{prop2}
Let $A\in\C^{n\times n}$ and consider the problem (\ref{epralhom}). Let $y(z)=(y_1(z),\ldots,y_n(z))$ be the solution of any Cauchy problem associated to equation (\ref{epralhom}). Then, 
\begin{itemize}
\item[-] If $\lambda=0$ is the only eigenvalue of $A$, then $y(z)$ has a polynomial growth at infinity. In addition to this, there exists $C>0$ such that
$$|y_j(z)|\le C|z|^{n-1},\qquad 1\le j\le n,\quad z\in\C.$$
\item[-] There exist $C_1,C_2>0$ such that 
$$|y_j(z)|\le C_1\exp(M(C_2|z|)),\qquad 1\le j\le n,\quad z\in\C,$$
where the function $M(\cdot)$ is defined in (\ref{e94}).
\end{itemize}
\end{prop}
\begin{proof}
The first part of the result is a direct consequence of the form of the explicit solution of (\ref{epralhom}), described in the proof of Theorem~\ref{teonodiag}. If $A\equiv 0$, then the statement is clear. Otherwise, one has that $y(z)$ is a linear combination of constant vectors multiplied by one of the elements in $\{E(\lambda z),\Delta_1E(\lambda z),\ldots,\Delta_{n-1}E(\lambda z)\}$, i.e. in $\{1/m(0),z/m(1),\ldots,z^{n-1}/m(n-1)\}$. The result follows directly from here.

For the second part of the proof we make use of Proposition~\ref{propk}, and observe that any solution of $y(z)$ is a linear combination of constant vectors multiplied by one of the elements in $\cup_{\lambda \in\hbox{spec}(A)}\{E(\lambda z),\Delta_1E(\lambda z),\ldots,\Delta_{h-1}E(\lambda z)\}$, for some $0\le h\le n-1$. We observe that, given any $\lambda\in\hbox{spec}(A)$, $0\le h\le n-1$, and $z\in\C$, one has that $\Delta_h E(\lambda z)=\sum_{p\ge0}a_pz^p$, with
$$|a_p|\le {p\choose h}\frac{|\lambda|^{p-h}}{m(p)}\le |\lambda|^{-h}(2|\lambda|)^{p}/m(p),$$
for some $h\ge0$. At this point, we recall the fact that the sequence $\mathbb{M}$ and the sequence of moments $m$ are equivalent (see Proposition~\ref{propequiv}), which entails that
$$|\Delta_h E(\lambda z)|\le C_{h,1}\exp(M(C_{h,2}|z|)),\qquad z\in\C,$$
for some $C_{h,1},C_{h,2}>0$. Therefore, $|y_j(z)|$ is upper bounded by a linear combination of functions of the form $\exp(M(C_{h}|z|))$, for certain positive constants $C_{h}$ and every $z\in\C$. The definition of $M(\cdot)$ allows us to group all of them into one term achieving  the statement in the result for some $C_1$ and where $C_2$ is the maximum of the constants $C_h$. 
\end{proof}

Further information of the growth of the solutions at infinity can be provided in terms of the measurement given by the order and type of the entire solutions. The following definition and results can be found in~\cite{maergoiz,maergoiz2}

\begin{defin}\label{def444}
Let $f\in\mathcal{O}(\C)$. We write $M_f(r):=\max\{|f(z)|: |z|=r\}$ for every $r\ge0$. The order of $f$ is defined by
$$\rho=\rho_f:=\displaystyle\limsup_{r\to\infty}\frac{\ln^{+}(\ln^{+}(M_{f}(r)))}{\ln(r)}.$$
Given $f\in\mathcal{O}(\C)$ of order $\rho\in\R$, the type of $f$ is defined by
$$\sigma=\sigma_f:=\displaystyle\limsup_{r\to\infty}\frac{\ln^{+}(M_f(r))}{r^{\rho}}.$$
Here, $\ln^+(\cdot)=\max\{0,\ln(\cdot)\}$.
\end{defin}

\begin{example}
The previous elements measure the growth of an entire function at infinity. In the case of $f(z)=\exp(\sigma z^{\rho})$, for some $\rho,\sigma>0$, one has that 
$$M_f(r)=\max\{\exp(\sigma r^{\rho}\cos(\rho\theta)):\theta\in\R\}= \exp(\sigma r^{\rho}),$$
and
$$\rho_f=\lim_{r\to\infty}\frac{\ln^+(\sigma r^{\rho})}{\ln(r)}=\rho,\qquad \sigma_f=\lim_{r\to\infty}\frac{\sigma r^{\rho}}{r^{\rho}}=\sigma.$$
\end{example}

\begin{lemma}[Theorem 4.2.1~\cite{holland}]\label{lema450}
Let $f_1,f_2\in\mathcal{O}(\C)$. Then $\rho_{f_1+f_2}=\max\{\rho_{f_1},\rho_{f_2}\}$.
\end{lemma}

\begin{lemma}[Theorem 4.2.3~\cite{holland}]\label{lema451}
Let $p\in\C[z]$ with $p\neq 0$, and $f\in\mathcal{O}(\C)$. Then, $\rho_{pf}=\rho_f$.
\end{lemma}

The following result is a straighforward application of the definition of type associated to a function.

\begin{lemma}\label{lema460}
Let $f,g\in\mathcal{O}(\C)$ both of order $\rho>0$. Then, $\sigma_{f+g}\le \max\{\sigma_f,\sigma_g\}$. For all $C\in\R^{\star}$, $\sigma_{Cf}=\sigma_f$.
\end{lemma}

The study of the order and type associated to an entire function of the form (\ref{e143}) is detailed in~\cite{maergoiz}, Section 3. More precisely, one has the following result whose writting has been adapted to our settings.

\begin{theo}[Theorem 3.3,~\cite{maergoiz}]\label{teo455}
Let $\mathbb{M}$ be a strongly regular sequence which admits a nonzero proximate order, say $\rho(t)\to\rho>0$, for $t\to\infty$. Then, the function $E(z)$ defined by (\ref{e143}) is an entire function of order $\rho$ and type 1.
\end{theo}

We also make use of the following result.

\begin{theo}[Theorem 1.10.3,~\cite{maergoiz2}]\label{teo456}
Let $f(z)=\sum_{p\ge0}a_pz^p$ be an entire function. Then,
\begin{itemize}
\item[i)] The order $\rho\ge 0$ of $f$ is given by 
$$\rho=\lim\sup_{p\to\infty}\frac{p\ln(p)}{-\ln|a_p|}.$$
\item[ii)] If $0<\rho<\infty$, then the type $\sigma\ge0$ of $f$ is determined by
$$(\sigma e \rho)^{1/\rho}=\limsup_{p\to\infty}p^{1/\rho}|a_p|^{1/p}.$$
\end{itemize}
\end{theo}

In view of Theorem~\ref{teo455} and Theorem~\ref{teo456}, one can give more information on the growth at infinity of the solutions to (\ref{epralhom}).

\begin{theo}\label{teo5}
Let $\mathbb{M}$ be a strongly regular sequence which admits a nonzero proximate order, say $\rho(t)\to\rho>0$, for $t\to\infty$. Let $y=y(z)$ be a solution of (\ref{epralhom}). Then, 
\begin{enumerate}
\item if $A$ admits a nonzero eigenvalue, then $y$ is an entire function of order $\rho$ and type upper bounded by $\sigma:=\max\{|\lambda|^{\rho}:\lambda\in\hbox{spec}(A)\}$, or an entire function of order 0.
\item if 0 is the only eigenvalue of $A$, then $y$ is a polynomial. Therefore its order is zero.
\end{enumerate}
\end{theo}
\begin{proof}
The second part of the statement is clear from the construction of the general solution of (\ref{epralhom}). Assume that $A$ admits a nonzero eigenvalue, say $\lambda$. First, we prove that for every $h\ge0$, the entire function $\Delta_hE(\lambda z)$ defined in (\ref{e258}) is of order $\rho$ and type $|\lambda|^{\rho}$, for every $\lambda\in\C^{\star}$. Let $h\ge0$ be an integer number and $\lambda\in\C^{\star}$. We observe that
\begin{equation}\label{e475}
\lim_{p\to\infty}\frac{-\ln\left[{p\choose h}|\lambda|^{p-h}\right]}{p\ln(p)}=\lim_{p\to\infty}\frac{-\ln(p^{h}|\lambda|^{p-h}h!)}{p\ln(p)}=0,
\end{equation}
which entails that
$$\limsup_{p\to\infty}\frac{p\ln(p)}{-\ln\left[{p\choose h}\frac{|\lambda|^{p-h}}{m(p)}\right]}=\limsup_{p\to\infty}\frac{p\ln(p)}{-\ln\left[\frac{1}{m(p)}\right]}=\rho,$$
regarding Theorem~\ref{teo455}. From Theorem~\ref{teo456} we conclude that $\Delta_hE(\lambda z)$ is an entire function of order $\rho$. On the other hand, we have that 
$$\lim_{p\to\infty}\left[{p\choose h}|\lambda|^{p-h}\right]^{1/p}=|\lambda|.$$
Therefore, in view of (\ref{teo455}) one has that
$$\limsup_{p\to\infty}p^{1/\rho}\left[{p\choose h}\frac{|\lambda|^{p-h}}{m(p)}\right]^{1/p}=\limsup_{p\to\infty}p^{1/\rho}\left[\frac{1}{m(p)}\right]^{1/p}|\lambda|=|\lambda|.$$
We conclude from Theorem~\ref{teo456} that the type of $\Delta_hE(\lambda z)$ is $|\lambda|^{\rho}$. Indeed, let $\sigma$ the type of the function $\Delta_hE(\lambda z)$. Then, one has that
$$(\sigma e\rho)^{1/\rho}=\limsup_{p\to\infty}p^{1/\rho}\left({p\choose h}\frac{|\lambda|^{p-h}}{m(p)}\right)^{1/p}=|\lambda|\limsup_{p\to\infty}p^{1/\rho}\left(\frac{1}{m(p)}\right)^{1/p}=|\lambda|(e\rho)^{1/\rho},$$
which entails that $\sigma=|\lambda|^{\rho}$.

In view of the form of the solution of (\ref{epralhom}) given by (\ref{e283}) and applying Lemma~\ref{lema450} and Lemma~\ref{lema451}, one gets that the order of $y$ is $\rho$. Observe that the order falls to zero if the coefficients associated to nonzero eigenvalues in the solution vanish. Lemma~\ref{lema460} allows to conclude that, in the case that the order is $\rho$, then the type of the solution is, at most $\sigma:=\max\{|\lambda|:\lambda\in\hbox{spec}(A)\}$ .
\end{proof}

\noindent\textbf{Remark:} Observe that the concrete type associated to a solution of a Cauchy problem associated to (\ref{epralhom}) is determined by the Cauchy data. For example, if $m=(p!)_{p\ge0}$ and $A=\hbox{diag}(1,2)$, the solution of the Cauchy problem with $y(0)=(1,0)$ has order 1 and type equal to 1, whereas the type of the solution of the Cauchy problem with $y(0)=(0,1)$ is of type 2. The same holds for the order. Indeed, consider the matrix $A=\hbox{diag}(1,0)$. The solution of the Cauchy problem with $y(0)=(0,1)$ is $y(z)\equiv(0,1)$ whereas the solution with $y(0)=(1,0)$ has order and type equal to 1.

\subsection{On the radial growth of the solutions}\label{sec51}

The radial growth properties of the kernel function $E$ allow us to give some information on the growth at infinity of the solutions to (\ref{epralhom}) along rays to infinity. More precisely, one has the following first result in this direction. 

\begin{prop}\label{prop4}
Let $\mathbb{M}$ be a strongly regular sequence which admits a nonzero proximate order $\rho(t)\to\rho>0$. Let $m$ be the sequence of moments constructed as in Definition~\ref{defi128}. Let $A\in\C^{n\times n}$ be a diagonalizable matrix and consider the Cauchy problem (\ref{e176}), for some Cauchy data $(x_0,y_0)\in\C^{1+n}$. Then, the solution $y=(y_1,\ldots,y_n)$ of (\ref{e176}) satisfies that for every $1\le h\le n$
$$|y_h(re^{i\theta})|\le \frac{C}{r^{\beta}},\qquad r\ge R_0$$
for some $C,\beta>0$ and some $R_0>0$ provided that $\theta\in\mathcal{A}_h$. $\mathcal{A}_h$ is some (possibly empty) set of directions.
\end{prop}
\begin{proof}
Let (\ref{e235}) be the general solution of (\ref{epralhom}). For any choice of Cauchy data $(z_0,y_0)\in\C^{1+n}$, let $C^0_{j,p}\in\C$ be the constants for $1\le j\le k$ and $1\le p\le \ell_j$ such that
\begin{equation}\label{e531}
y_h(z)=\sum_{j=1}^{k}\sum_{p=1}^{\ell_j}C^{0}_{j,p}E(\lambda_j z)v_{j,p,h},\quad 1\le h\le n
\end{equation}
determines the unique solution to the Cauchy problem (\ref{e176}). The set $\mathcal{A}_h$ is defined by
$$\mathcal{A}_{h}=\bigcap_{j=1}^{k}\left\{\theta\in\R:\frac{\omega(\mathbb{M})\pi}{2}-\arg(\lambda_j)<\theta<2\pi-\frac{\omega(\mathbb{M})\pi}{2}-\arg(\lambda_j)\right\}.$$
Observe that $\theta\in\mathcal{A}_h$ if and only if $re^{i\theta}\lambda_j\in S_{\pi}(\nu)$, for some $0<\nu<2\pi-\omega(\mathbb{M})\pi$, for every $1\le j\le k$ which entails bounds on the kernel function as in (\ref{e145}).
\end{proof}

The previous result distinguishes certain directions for which the bound established in Proposition~\ref{prop2} has been tightened. The following definition and results can be found in~\cite{maergoiz,maergoiz2}. They will be the key point to determine the growth of the solution to (\ref{epralhom}) along rays to infinity.

\begin{defin}
Let $\rho(t)$ be a proximate order of $f\in\mathcal{O}(\C)$, i.e. $\rho(t)$ is a proximate order and 
$$0<\sigma_f:=\limsup_{r\to\infty}\frac{\ln(M_f(r))}{r^{\rho(r)}}<\infty.$$
Then, $\sigma_f$ is known as the type of $f$ relative to $\rho$, and
$$h_f(\theta)=\limsup_{r\to\infty}\frac{\ln|f(re^{i\theta})|}{r^{\rho(r)}},\quad\theta\in\R,$$
is known as the generalized indicator of $f$. We write $h_{f}^{+}=\max\{0,h_f\}$.
\end{defin}

\begin{theo}[Theorem 3.7,~\cite{maergoiz}]
Let $\mathbb{M}$ be a strongly regular sequence which admits a proximate order $\rho(t)\to\rho>0$. Let $E$ be the associated kernel function defined in (\ref{e143}). Then,
\begin{itemize}
\item[i)] 
$$h_E^+(\theta)=\left\{ \begin{array}{ccc}
             \cos(\rho\theta) &   if\quad &  |\theta|\le \min\{\pi,\frac{\pi}{2\rho}\}\\
             0  &  if \quad &  \rho\ge \frac{1}{2},\quad \min\{\pi,\frac{\pi}{2\rho}\}\le |\theta|\le \pi
             \end{array}
   \right.$$
for every $\theta\in[-\pi,\pi]$.
\item [ii)] If $\rho>1$, then $h_E(\theta)=h^+_E(\theta)$, for $\theta\in\R$.
\item[iii)] If $\rho\le 1/2$, then $h_E^+(\theta)=\cos(\rho\theta)$, $|\theta|\le \pi$.
\item[iv)] If $1/2<\rho\le 1$ and there exists $k\in\N$ such that $m(-k-1)\neq0$, then $h_E(\theta)=h_E^+(\theta)$, $\theta\in\R$.
\end{itemize}
\end{theo}


\begin{lemma}\label{lema559}
Let $\rho(t)$ be a proximate order of $f\in\mathcal{O}(\C)$. Then, for every $C\in\C^{\star}$ the function $t\mapsto\rho(|C|t)$ is a proximate order of the function $t\mapsto f(Ct)$. Moreover, the type of the function $t\mapsto f(Ct)$ is given by $\sigma_{f}|C|^{\rho}$.
\end{lemma}
\begin{proof}
It is straight to check that the function $\tilde{\rho}(t)=\rho(|C| t)$ is a proximate order. Observe that (\ref{e159}) holds because
$$\lim_{r\to\infty}\tilde{\rho}(r)=\lim_{r\to\infty}\rho(r)=\rho.$$
In addition to this, one has
$$\lim_{r\to\infty}r\tilde{\rho}'(r)\ln(r)=\lim_{r\to\infty}|C|r\rho'(|C|r)\ln(|C|r)\frac{\ln(r)}{\ln(|C|r)}.$$
Due to $\frac{\ln(r)}{\ln(|C|r)}\to 1$ for $r\to\infty$, we get that the previous limit is zero because $\rho(t)$ is a proximate order. This entails that (\ref{e169}) holds and $\tilde{\rho}$ is a proximate order. 

We now prove that $\tilde{\rho}(t)$ is a proximate order of $\tilde{f}(t)=f(Ct)$. First, observe that
\begin{align*}
M_{f(Cz)}(r)&=\sup_{|z|=r}\{|f(\lambda z)|\}=\sup_{\theta\in\R}\{|f(|\lambda|re^{i(\theta+\arg(\lambda))})\}\\
&=\sup_{\theta\in\R}\{|f(|\lambda|re^{i\theta})\}=\sup_{|z|=|\lambda| r}\{|f(z)|\}=M_{f}(|\lambda|r).
\end{align*}

From Definition~\ref{def444}, we have
\begin{equation}\label{e578}\limsup_{r\to\infty}\frac{\ln(M_{\tilde{f}}(r))}{r^{\tilde{\rho}(r)}}=\limsup_{r\to\infty}\frac{\ln(M_{f}(|C|r))}{r^{\rho(|C|r)}}=\limsup_{r\to\infty}\frac{\ln(M_{f}(|C|r))}{(|C|r)^{\rho(|C|r)}}|C|^{\rho(|C|r)}.
\end{equation}

From the properties of proximate order, one has $|C|^{\rho(|C|r)}\to |C|^{\rho}$ for $t\to\infty$. Therefore, the (\ref{e578}) equals $\sigma_{f}|C|^{\rho}>0$.
\end{proof}

\begin{lemma}\label{lema591}
Let $\lambda\in\C^{\star}$. Let $\mathbb{M}$ be a strongly regular sequence admitting a nonzero proximate order, say $\rho(t)\to \rho>0$, for $t\to\infty$. We consider the kernel associated kernel  function $E$ described in (\ref{e143}). Then, the generalized indicator of the function $z\mapsto E(\lambda z)$ is given by
$$h_{E(\lambda z)}(\theta)=|\lambda|^{\rho}h_{E}(\theta+\arg(\lambda)).$$
\end{lemma}
\begin{proof}
In view of Lemma~\ref{lema559}, the function $\rho(|\lambda|t)$ is a proximate order of the entire function $E(\lambda z)$.  We have
\begin{multline}
h_{E(\lambda z)}(\theta)=\limsup_{r\to\infty}\frac{\ln|E(\lambda r e^{i\theta})|}{r^{\rho(|\lambda|r)}}=\limsup_{r\to\infty}\frac{\ln|E(\lambda r e^{i\theta})|}{(|\lambda|r)^{\rho(|\lambda|r)}}|\lambda|^{\rho(|\lambda|r)}\\
=|\lambda|^{\rho}\limsup_{r\to\infty}\frac{\ln|E(r e^{i(\theta+\arg(\lambda)})|}{r^{\rho(r)}}=|\lambda|^{\rho}h_{E}(\theta+\arg(\lambda)).
\end{multline}
This concludes the proof.
\end{proof}

As a direct application of the previous Lemma, we come up with the explicit form of the generalized indicator of $E(\lambda z)$.

\begin{corol}
Let $\mathbb{M}$ be a strongly regular sequence which admits a proximate order $\rho(t)\to\rho>0$. Let $E$ be the associated kernel function defined in (\ref{e143}), and $\lambda\in\C^{\star}$. Then,
\begin{itemize}
\item[i)] 
\begin{equation}\label{e608}
h_{E(\lambda z)}^+(\theta)=\left\{ \begin{array}{ccc}
             |\lambda|^{\rho}\cos(\rho(\theta+\arg(\lambda))) &   if\quad &  |\theta+\arg(\lambda)|\le \min\{\pi,\frac{\pi}{2\rho}\}\\
             0  &  if \quad &  \rho\ge \frac{1}{2},\quad \min\{\pi,\frac{\pi}{2\rho}\}\le |\theta+\arg(\lambda)|\le \pi
             \end{array}
   \right.
	\end{equation}
for every $\theta\in[-\pi-\arg(\lambda),\pi-\arg(\lambda)]$.
\item [ii)] If $\rho>1$, then $h_{E(\lambda z)}(\theta)=h^+_{E(\lambda z)}(\theta)$, for $\theta\in\R$.
\item[iii)] If $\rho\le 1/2$, then $h_{E(\lambda z)}^+(\theta)=|\lambda|^{\rho}\cos(\rho(\theta+\arg(\lambda))$, $|\theta+\arg(\lambda)|\le \pi$.
\item[iv)] If $1/2<\rho\le 1$ and there exists $k\in\N$ such that $m(-k-1)\neq0$, then $h_{E(\lambda z)}(\theta)=h_{E(\lambda z)}^+(\theta)$, $\theta\in\R$.
\end{itemize}
\end{corol}
\begin{proof}
We observe that $h^{+}_{E(\lambda z)}(\theta)$ for every $\theta\in[-\pi-\arg(\lambda),\pi-\arg(\lambda)]$ is determined by
\begin{align*}
h^{+}_{E(\lambda z)}(\theta)&=\max\{h_{E(\lambda z)}(\theta),0\}\\
&=|\lambda|^{\rho}\max\{h_E(\theta+\arg(\lambda)),0\}=|\lambda|^{\rho}h^{+}_{E}(\theta+\arg(\lambda)),
\end{align*}
which yields (\ref{e608}). 

If $\rho>1$, then for every $\theta\in\R$ 
$$h_{E(\lambda z)}(\theta)=|\lambda|^{\rho}h_{E}(\theta+\arg(\lambda))=|\lambda|^{\rho}h^+_E(\theta+\arg(\lambda))=\max\{0,|\lambda|^{\rho}h_E(\theta+\arg(\lambda))\}=h^+_{E(\lambda z}(\theta).$$

If $\rho\le 1/2$, then for every $|\theta+\arg(\lambda)|\le \pi$
$$h_{E(\lambda z)}^+(\theta)=|\lambda|^{\rho}h^{+}_{E}(\theta+\arg(\lambda))=|\lambda|^{\rho}\cos(\rho(\theta+\arg(\lambda)).$$

Finally, if $1/2<\rho\le 1$ and there exists $k\in\N$ such that $m(-k-1)\neq 0$, then for every $\theta\in\R$ one has
$$h_{E(\lambda z)}(\theta)=|\lambda|^{\rho}h_{E}(\theta+\arg(\lambda))=|\lambda|^{\rho}h^{+}_{E}(\theta+\arg(\lambda))=h^+_{E(\lambda z)}(\theta).$$
\end{proof}

Our final result in this work states an upper bound for the generalized indicator of the solutions to (\ref{epralhom}).

\begin{theo}\label{teo6}
Let $\mathbb{M}$ be a strongly regular sequence which admits a nonzero proximate order $\rho(t)\to\rho>0$. Let $m$ be the sequence of moments constructed as in Definition~\ref{defi128}. Assume that $A\in \C^{n\times n}$ is a diagonalizable matrix, and consider the Cauchy problem (\ref{e176}), for some Cauchy data $(x_0,y_0)\in\C^{1+n}$. Then, the solution $y=(y_1,\ldots,y_n)$ of (\ref{e176}) satisfies for every $1\le h\le n$ that
$$h_{y_{h}}(\theta)\le \max\left\{|\lambda|^{\rho}h_{E}(\theta+\arg(\lambda)): \lambda\in\hbox{spec}(A)\right\},$$
for every $\theta\in\R$.

\end{theo}
\begin{proof}

Let $1\le h\le n$. We write $y_h$ in the form (\ref{e531}), for certain $C_{j,p}^0,v_{j,p,h}\in\C$, that we assume to be nonzero without loss of generality, and where $\{\lambda_{j}:1\le j\le k\}$ is the set of eigenvalues of $A$. We can also assume that $\lambda=0$ is not the only eigenvalue of $A$. Otherwise, $A\equiv 0$.

At this point, we may assume that at least one of the constants $C_{j,p}^0$ associated to a nonzero eigenvalue differs from zero. Otherwise, the solution of the Cauchy problem is a polynomial.

 We observe that for every $1\le j\le k$ the function $C^{0}_{j,p}E(\lambda_j z)v_{j,p,h}$ admits $t\mapsto \rho(|\lambda_j| t)$ as proximate order, with $\rho(t)$ being a proximate order associated to the kernel function $E$. In addition to this, the type of $C^{0}_{j,p}E(\lambda_j z)v_{j,p,h}$ coincides with the type of $E(\lambda_j z)$, which equals $|\lambda_j|^{\rho}$. 

For every $r\ge0$, one has
\begin{multline}
M_{y_{h}}(r)=\sup\left\{\left|\sum_{j=1}^{k}\sum_{p=1}^{\ell_j}C^{0}_{j,p}E(\lambda_j z)v_{j,p,h}\right|:|z|=r\right\}\\
\le \left( \sum_{j=1}^{k}\sum_{p=1}^{\ell_j}    |C^{0}_{j,p}v_{j,p,h}|\right)\sup\{\max_{1\le j\le k}|E(\lambda_j z)|:|z|=r\}.
\end{multline}
Let us write $C_1:=\sum_{j=1}^{k}\sum_{p=1}^{\ell_j}    |C^{0}_{j,p}v_{j,p,h}|$. We have obtained that 
$$M_{y_h}(r)\le C_1\max_{1\le j\le k}\{M_{E(\lambda_j z)}(r)\}.$$ 

Now, let $1\le j_0\le k$ be such that $|\lambda_{j_0}|\ge |\lambda_j|$ for every $1\le j\le k$.  We recall that $M_{E(\lambda_j z)}(r)=M_E(|\lambda_j|r)$ for every $r\ge0$. The monotonicity of $M_{E}$ entails that $M_{y_h}(r)\le C_1M_{E(\lambda_{j_0} z)}(r)$. We recall from the proof of Lemma~\ref{lema559} that $t\mapsto \rho(|\lambda_{j_0}|t)$ is a proximate order of $E(\lambda_{j_0}z)$. In order to check that it is also a proximate order for $y_{h}$ we have
\begin{multline}
\limsup_{r\to\infty}\frac{\ln(M_{y_{h}}(r))}{r^{\rho(|\lambda_0|r)}}\le \limsup_{r\to\infty}\frac{\ln(C_1M_{E}(|\lambda_{j_0}|r))}{r^{\rho(|\lambda_0|r)}}
=\limsup_{r\to\infty}\frac{\ln(M_{E}(|\lambda_{j_0}|r))}{r^{\rho(|\lambda_0|r)}}\\
=\limsup_{r\to\infty}\frac{\ln(M_{E}(|\lambda_{j_0}|r))}{|\lambda_0|r)^{\rho(|\lambda_0|r)}}|\lambda_0|^{\rho(|\lambda_0|r)}=\sigma_{E}|\lambda_0|^{\rho}=|\lambda_0|^{\rho}.
\end{multline}

Therefore, $\sigma_{y_h}\le |\lambda_0|^{\rho}$.

Now, we study 
$$\limsup_{r\to\infty}\frac{\ln|y_h(re^{i\theta})|}{r^{\rho(|\lambda_0|r)}}.$$
Analogous computations as in the first part of the proof, and taking into account the properties of $\limsup$ yield to upper estimate the previous expression by
$$\limsup_{r\to\infty}\frac{\ln(C_1n)+\max_{1\le j\le k}\ln|E(\lambda_j r e^{i\theta})|}{r^{\rho(|\lambda_0|r)}}\le\max_{1\le j\le k}\left\{\limsup_{r\to\infty}\frac{\ln|E(\lambda_j r e^{i\theta})|}{r^{\rho(|\lambda_0|r)}}\right\}.$$


taking into account Lemma~\ref{lema591}we arrive at
\begin{multline}
\limsup_{r\to\infty}\frac{\ln|E(\lambda_j r e^{i\theta})|}{r^{\rho(|\lambda_0|r)}}
=\limsup_{r\to\infty}|\lambda_0|^{\rho(|\lambda_0|r)}\frac{\ln|E(\frac{|\lambda_j|}{|\lambda_0|}|\lambda_0| r e^{i(\theta+\arg(\lambda_j))})|}{(|\lambda_0|r)^{\rho(|\lambda_0|r)}}\\
=|\lambda_0|^{\rho}h_{E(|\lambda_j|/|\lambda_0|z)}(\theta+\arg(\lambda_j))=|\lambda_j|^{\rho}h_{E}(\theta+\arg( \lambda_j).
\end{multline}

This entails that
$$\limsup_{r\to\infty}\frac{\ln|y_h(re^{i\theta})|}{r^{\rho(|\lambda_0|r)}}\le \max_{1\le j\le k}|\lambda_j|^{\rho}h_{E}(\theta+\arg( \lambda_j),$$
arriving at the conclusion.
\end{proof}

\noindent\textbf{Remark:} Observe that in the previous proof, we have also obtained a proximate order of the solution $y_h$ for all $1\le h\le n$ and also an upper bound for its type $\sigma_{y_h}$. Indeed, if $\lambda$ is the eigenvalue of $A$ of larger modulus, then the function $\rho(|\lambda|t)$ is a proximate order of $y_h$, and the type of $y_h$ is upper bounded by $|\lambda|^{\rho}$.

\vspace{0.3cm}

\noindent\textbf{Remark:} The exact type and order associated to a solution of a Cauchy problem are determined by the Cauchy data.






\begin{thebibliography}{99}
\bibitem{aghror} R. Agarwal, S. Hristova, D O'Regan, \emph{Mittag-Leffler stability for impulsive Caputo fractional differential equations}, Differ. Equ. Dyn. Syst. 29, No. 3 (2021) 689--705. 
\bibitem{ba2} W. Balser, \emph{Formal power series and linear systems of meromorphic ordinary differential equations.} Universitext. Springer-Verlag, New York, 2000. xviii+299 pp.
\bibitem{bayo} W. Balser, M. Yoshino, \emph{Gevrey order of formal power series solutions of inhomogeneous partial differential equations with constant coefficients.} Funkcial. Ekvac. 53 (2010) 411--434.
\bibitem{bonillaetal} B. Bonilla, M. Rivero, J. J. Trujillo, \emph{On systems of linear fractional differential equations with constant coefficients.} Appl. Math. Comput. 187 (2007) No. 1, 68--78. 
\bibitem{goos} A. A. Goldberg, I. V. Ostrovskii, Value Distribution of Meromorphic Functions, Translation of Mathematical Monograph, vol. 236. American Mathematical Society, Providence, 2008.
\bibitem{holland} A. S. B. Holland, Introduction to the theory of entire functions. Pure and Applied Mathematics, 56. New York-London: Academic Press, 1973. 
\bibitem{immink1} G. K. Immink, \emph{Exact asymptotics of nonlinear difference equations with levels 1 and $1+$}, Ann. Fac. Sci. Toulouse 2 (2008), 309--356.
\bibitem{immink2} G. K. Immink, \emph{Accelero-summation of the formal solutions of nonlinear difference equations}, Ann. Inst. Fourier (Grenoble) 61(1) (2011) 1--51.
\bibitem{jikalasa} J. Jim\'enez-Garrido, S. Kamimoto, A. Lastra, J. Sanz, \emph{Multisummability in Carleman ultraholomorphic classes by means of nonzero proximate orders,} J. Math. Anal. Appl. 472, No. 1 (2019) 627--686. 
\bibitem{jisash} J. Jim\'enez-Garrido, J. Sanz, G. Schindl, \emph{Injectivity and surjectivity of the asymptotic Borel map in Carleman ultraholomorphic classes.} J. Math. Anal. Appl. 469 (2019), 136--168.
\bibitem{jisash2} J. Jim\'enez-Garrido, J. Sanz, G. Schindl, \emph{Log-convex sequences and nonzero proximate orders.} J. Math. Anal. Appl. 448(2) (2017) 1572--1599. 
\bibitem{komatsu} h. Komatsu, \emph{Ultradistributions. I: Structure theorems and a characterization}, J. Fac. Sci., Univ. Tokyo, Sect. I A 20 (1973) 25--105. 
\bibitem{lamasa} A. Lastra, S. Malek, J. Sanz, \emph{Summability in general Carleman ultraholomorphic classes.} J. Math. Anal. Appl. 430 (2015) 1175--1206. 
\bibitem{lamasa2} A. Lastra, S. Malek, J. Sanz, \emph{Strongly regular multi-level solutions of singularly perturbed linear partial differential equations,} Results Math. 70 (2016), no. 3--4, 581--614.
\bibitem{lamisu} A. Lastra, S. Michalik, M. Suwi\'nska, \emph{Summability of formal solutions for some generalized moment partial differential equations.} Result. Math. 76, No. 1 (2021) Paper No. 22.
\bibitem{lamisu2} A. Lastra, S. Michalik, M. Suwi\'nska, \emph{Estimates of formal solutions for some generalized moment partial differential equations,} J. Math. Anal. Appl. 500 (2021), no. 1.
\bibitem{lamisu3} A. Lastra, S. Michalik, M. Suwi\'nska, \emph{Summability of formal solutions for a family of generalized moment integro-differential equations,} to appear in Fract. Calc. Appl. Anal., 2021.
\bibitem{maergoiz} L. S. Maergojz, \emph{Indicator diagram and generalized Borel-Laplace transforms for entire functions of a given proximate order}, St. Petersbg. Math. J. 12, No. 2 (2001) 191--232; translation from Algebra Anal. 12, No. 2  (2000) 1--63. 
\bibitem{maergoiz2} L. S. Maergojz, Asymptotic characteristics of entire functions and their applications in mathematics and biophysics. Mathematics and its Applications (Dordrecht) 559. Dordrecht: Kluwer Academic Publishers, xxiv, 2003. 
\bibitem{mandel} S. Mandelbrojt, \emph{S\'eries Adh\'erentes, Regularisation des suites, Applications}, Gauthier-Villars, Paris, 1952.
\bibitem{maab} M. M. Matar, E. S. A. Skhail, \emph{On stability analysis of semi-linear fractional differential systems,} Math. Methods Appl. Sci. 43, No. 5 (2020) 2528--2537. 
\bibitem{matignon} D. Matignon, \emph{Stability results for fractional differential equations with applications to control processing,} Proc. Comput. Eng. Syst. Appl. 2 (1996), 963--968.Computational Engineering in System Application, Vol. 2 1996, p. 963. 
\bibitem{maon} I. Matychyn, V. Onyshchenko, \emph{Optimal control of linear systems with fractional derivatives,} Fract. Calc. Appl. Anal. 21, No. 1 (2018) 134--150. 
\bibitem{mi} S. Michalik, \emph{Analytic solutions of moment partial differential equations with constant coefficients,} Funkcial. Ekvac. 56 (2013), no. 1, 19--50.
\bibitem{michalik12} S. Michalik, \emph{Multisummability of formal solutions of inhomogeneous linear partial differential equations with constant coefficients}, J. Dyn. Control Syst. 18 (2012) 103--133.
\bibitem{mitk} S. Michalik, B. Tkacz, \emph{The Stokes phenomenon for some moment partial differential equations,} J. Dyn. Control Syst. 25 (2019), no. 4, 573--598.
\bibitem{sanz} J. Sanz, \emph{Flat functions in Carleman ultraholomorphic classes via proximate orders}, J. Math. Anal. Appl. 415(2) (2014), 623--643.
\bibitem{sanzrev} J. Sanz, \emph{Asymptotic analysis and summability of formal power series,} Analytic, algebraic and geometric aspects of differential equations, 199--262, Trends Math., Birkh\"auser/Springer, Cham, 2017.
\bibitem{shafarevich} I. R. Shafarevich, A. O. Remizov, Linear Algebra and Geometry, Springer, 2012.
\bibitem{su} M. Suwi\'nska, \emph{Gevrey estimates of formal solutions for certain moment partial differential equations with variable coefficients,} J. Dyn. Control Syst. 27, No. 2 (2021) 355--370. 
\bibitem{thilliez} V. Thilliez, \emph{Division by flat ultradifferentiable functions and sectorial extensions}, Result. Math. 44 (2003), 169--188.
\bibitem{ucoz} E. U\c{c}ar, N. \"Ozdemir, \emph{A fractional model of cancer-immune system with Caputo and Caputo–Fabrizio derivatives,} Eur. Phys. J. Plus 136 (2021) 43.
\end{thebibliography}
\end{document}